\documentclass[4pt, oneside]{article}
\usepackage{amsmath,amssymb,amsthm,units,stmaryrd}

\usepackage{yfonts}
\usepackage{units,stmaryrd}
\usepackage{color}
\usepackage{graphicx}
\usepackage[all]{xy}

\newtheorem{definition}{Definition}[section]
\newtheorem{lemm}{Lemma}[section]

\newtheorem{prop}{Proposition}[section]
\newtheorem{cor}{Corollary}[section]

\newtheorem{theorem}{Theorem}[section]

\def\lb{\left\llbracket}
\def\rb{\right\rrbracket}

\newcommand{\glp}{{\ensuremath{\mathsf{GLP}}}\xspace}

\usepackage{xspace}

\def\le{{\ell}}
\def\ex{e}

\def\lb{\left\llbracket}
\def\rb{\right\rrbracket}

\def\fmodels{\xymatrix{
\ar@{|=}[r]^{<\omega}&
}
}
\def\nmodels{\xymatrix{
\ar@{|=}[r]^{N}&
}
}
\def\<{\left <}

\def\nc{{\Box}}
\def\ps{{\Diamond}}

\def\>{\right >}
\def\bra {\left [}
\def\ket{\right ]}
\def\cbra{\left \{}
\def\cket{\right \}}

\DeclareSymbolFont{AMSb}{U}{msb}{m}{n}
\DeclareMathSymbol{\N}{\mathbin}{AMSb}{"4E}
\DeclareMathSymbol{\Z}{\mathbin}{AMSb}{"5A}
\DeclareMathSymbol{\R}{\mathbin}{AMSb}{"52}
\DeclareMathSymbol{\Q}{\mathbin}{AMSb}{"51}
\DeclareMathSymbol{\I}{\mathbin}{AMSb}{"49}
\DeclareMathSymbol{\C}{\mathbin}{AMSb}{"43}

\newcommand{\comment}[2]{{\color{red}{\bf #2}}}

\begin{document}

\title{Models of transfinite provability logic
      }

\author{David Fern\'{a}ndez-Duque\footnote{Group for Computational Logic, Universidad de Sevilla, dfduque@us.es} \and Joost J. Joosten\footnote{Department of Logic, History and Philosophy of Science, University of Barcelona, jjoosten@ub.edu}}

\maketitle

\begin{abstract}
For any ordinal $\Lambda$, we can define a polymodal logic $\mathsf{GLP}_\Lambda$, with a modality $[\xi]$ for each $\xi<\Lambda$. These represent provability predicates of increasing strength. Although $\mathsf{GLP}_\Lambda$ has no Kripke models, Ignatiev showed that indeed one can construct a Kripke model of the variable-free fragment with natural number modalities, denoted $\mathsf{GLP}^0_\omega$. Later, Icard defined a topological model for $\mathsf{GLP}^0_\omega$ which is very closely related to Ignatiev's.

In this paper we show how to extend these constructions for arbitrary $\Lambda$. More generally, for each $\Theta,\Lambda$ we build a Kripke model $\mathfrak I^\Theta_\Lambda$ and a topological model $\mathfrak T^\Theta_\Lambda$, and show that $\mathsf{GLP}^0_\Lambda$ is sound for both of these structures, as well as complete, provided $\Theta$ is large enough.
\end{abstract}


\section{Introduction}
It was G\"odel who first suggested interpreting the modal $\nc$ as a provability predicate, which as he observed should satisfy $\nc(\phi\to\psi)\to(\nc\phi\to\nc\psi)$ and $\nc\phi\to\nc\nc\phi$. With this, the Second Incompleteness Theorem could be expressed succinctly as $\ps\top\to\ps\nc\bot$.

More generally, L\"ob's axiom $\nc(\nc\phi\to\phi)\to\nc\phi$ is valid for this interpretation, and with this we obtain a complete characterization of the propositional behavior of provability in Peano Arithmetic \cite{Solovay:1976}. The modal logic obtained {from L\"ob's axiom} is called $\mathsf{GL}$ (for G\"odel-L\"ob) and is rather well-behaved; it is decidable and has finite Kripke models, based on transitive, well-founded frames \cite{Segerberg:1971}.

Japaridze \cite{Japaridze:1988} then suggested extending $\mathsf{GL}$ by a sequence of provability modalities $[n]$, for $n<\omega$, where {$[n]\phi$} could be interpreted (for example) as {\em $\phi$ is derivable using $n$ instances of the $\omega$-rule}. We shall refer to this extension as $\mathsf{GLP}_\omega$. $\mathsf{GLP}_\omega$ turns out to be much more powerful than $\mathsf{GL}$, and indeed Beklemishev has shown how it can be used to perform ordinal analysis of Peano Arithmetic and its natural subtheories \cite{Beklemishev:2004}.

However, as a modal logic, it is much more ill-behaved than $\mathsf{GL}$. Most notably, over the class of $\mathsf{GLP}$ Kripke frames, the formula $[1]\bot$ is valid! This is clearly undesirable. There are ways to get around this, for example using topological semantics. However, Ignatiev {in \cite{Ignatiev:1993}} showed how one can still get Kripke frames for the {\em closed} fragment of $\mathsf{GLP}_\omega$, which contains no propositional variables (only $\bot$). This fragment, which we denote $\mathsf{GLP}^0_\omega$, is still expressive enough to perform Beklemishev's ordinal analysis.

Later, Icard provided topological models for $\mathsf{GLP}^0_\omega$ \cite{Icard:2009}. {T}he full logic actually does have topological models, and indeed has been proven complete for these semantics by Beklemishev and Gabelaia \cite{BeklGabel:2011}. However, this requires rather heavy machinery and some non-constructive {methods}, all of which can be avoided when dealing only with the closed fragment.

Our goal is to extend {the} results {on $\mathsf{GLP}^0_\omega$} to $\mathsf{GLP}^0_\Lambda$, where $\Lambda$ is an arbitrary ordinal (or, if one wishes, the class of all ordinals). To do this we build upon known techniques, but dealing with transfinite modalities poses many new challenges. In particular, models will now have to be much `deeper' if we wish to obtain completeness.

The layout of the paper is as follows. In Section \ref{syntax} we give a quick overview of the logics $\mathsf{GLP}^0_\Lambda$. Section \ref{motivation} then gives some motivation for the constructions we shall present.

In Section \ref{super} we discuss how one `hyperates' ordinal exponentiation and last exponents. {Hyperations are a form of transfinite iteration and} will be crucial in describing our models.

In Section \ref{les} we introduce {\em $\le$-sequences}, which provide a generalization of the ``worlds'' in the Kripke semantics of $\mathsf{GLP}^0_\omega$ introduced by Ignatiev. Then, Section \ref{gim} defines generalizations of Ignatiev models with arbitrary ``depth'' and ``length'' and shows that indeed they provide models for $\glp^0_\Lambda$.

In Section \ref{topological} we define topological models for $\mathsf{GLP}^0_\Lambda$; these are generalizations of the polytopological spaces introduced by Icard. Finally, Section \ref{seccomp} proves soundness and establishes conditions on these models under which $\glp^0_\Lambda$ is complete for them.

\section{The logic $\mathsf{GLP}^0_\Lambda$}\label{syntax}

Let $\Lambda$ be either an ordinal or the class of all ordinals. Formulas of $\mathsf{GLP}^0_\Lambda$ are built from $\bot$ using Boolean connectives $\neg,\wedge$ and a modality $[\xi]$ for each $\xi<\Lambda$. As is customary, we use $\<\xi\>$ as a shorthand for $\neg[\xi]\neg$.

Note that there are no propositional variables, as we are concerned here with the {\em closed fragment} of $\mathsf{GLP}_\Lambda$.

The logic $\mathsf{GLP}^0_\Lambda$ is given by the following rules and axioms:
\begin{enumerate}
\item all propositional tautologies{,}
\item $[\xi](\phi\to\psi)\to([\xi]\phi\to[\xi]\psi)$ for all $\xi<\Lambda${,}
\item {$[\xi]([\xi]\phi \to \phi)\to[\xi]\phi$ for all $\xi<\Lambda$}{,}
\item $\<\zeta\>\phi\to\<\xi\>\phi$ for $\xi<\zeta<\Lambda${,}
\item $\<\xi\>\phi\to [\zeta]\<\xi\>\phi$ for $\xi<\zeta<\Lambda$.
\end{enumerate}

A {\em Kripke frame} is a structure $\mathfrak F=\<W,\<R_i\>_{i<I}\>$, where $W$ is a set and $\<R_i\>_{i<I}$ a family of binary relations on $W$. Since we are restricting to the closed fragment we make no distinction between Kripke frames and Kripke {\em models}. To each formula $\psi$ in the closed modal language with modalities $\<i\>$ for $i<I$ we assign a set $\lb\psi\rb_\mathfrak F\subseteq W$ inductively as follows:
\[
\begin{array}{lcl}
\lb\bot\rb_\mathfrak F&=&\varnothing\\\\
\lb\neg\phi\rb_\mathfrak F&=&W\setminus\lb\phi\rb_\mathfrak F\\\\
\lb\phi\wedge\psi\rb_\mathfrak F&=&\comment{W}{}\lb\phi\rb_\mathfrak F\cap\lb\psi\rb_\mathfrak F\\\\
\lb\<i\>\phi\rb_\mathfrak F&=&R^{-1}_i\lb\phi\rb_\mathfrak F.
\end{array}
\]

Often we will write $\<\mathfrak F,x\>\models\psi$ instead of $x\in\lb \psi\rb_\mathfrak F$.

It is well-known that polymodal $\mathsf{GL}$ is sound for $\mathfrak F$ whenever $R^{-1}_i$ is well-founded and transitive, in which case we write it $<_i$. However, constructing models of $\glp_\Lambda$ is substantially more difficult than constructing models of $\mathsf{GL}$, as we shall see. 


\section{Motivation for our model}\label{motivation}


The full logic $\mathsf{GLP}_\Lambda$ cannot be sound and complete with respect to any class of Kripke frames. Indeed, let $\mathfrak F=\langle W,\<<_\xi\>_{\xi<\lambda}\rangle$ be a polymodal frame.

Then, it is not too hard to check that
\begin{enumerate}
\item L\"ob's axiom $[\xi]([\xi]\phi \to \phi)\to [\xi]\phi$ is valid if and only if $<_\xi$ is well-founded and transitive{,}

\item the axiom $\langle {\zeta} \rangle \phi \to \langle {\xi} \rangle \phi$ {for $\xi\leq\zeta$} is valid if and only if, whenever $w<_{\zeta}v${}, then $w<_{\xi}v$, and

\item $\langle \xi \rangle \phi \to [\zeta] \langle \xi \rangle \phi$ {for $\xi<\zeta$} is valid if, whenever $v<_{\zeta}w $, $u<_{\xi}w$ and $\xi<\zeta$, then $u<_{\xi}v$.
\end{enumerate}
Suppose that for $\xi<\zeta$, there are two worlds such that $w<_\zeta v$. Then from 2 we see that $w<_\xi v$, while from 3 this implies that $w<_\xi w$. But this clearly violates 1. Hence if $\mathfrak F\models \mathsf{GLP}$, it follows that all accessibility relations (except possibly $<_0$) are empty.

However, this does not rule out the possibility that the closed fragments $\mathsf{GLP}^0_\Lambda$ have Kripke frames for which they are sound and complete. This turned out to be the case for $\mathsf{GLP}^0_{\omega}$ and in the current paper we shall extend this result to $\mathsf{GLP}^0_\Lambda$, with $\Lambda$ arbitrary.

More precisely, given ordinals $\Lambda,\Theta$, we will construct a Kripke frame $\mathfrak I^\Theta_\Lambda$ with `depth' $\Theta$ (i.e., the order-type of $<_0$) and `length' $\Lambda$ (the set of modalities it interprets). $\mathfrak I^\Theta_\Lambda$ validates all frame conditions except for condition 3. We shall only approximate it in that we require, for $\xi<\zeta$,
\[
v<_{\zeta}w \Rightarrow \exists\, v'<_\xi w \text{ such that }v'\leftrightarroweq_{\mathsf p} v.
\]
Here $\mathsf p$ will be a set of parameters and $u'\leftrightarroweq_{\mathsf p}u$ denotes that $u'$ is {\em $\mathsf p$-bisimilar} to $u$. The parameters $\mathsf p$ can be adjusted depending on $\phi$ in order to validate each instance of the axiom.



One convenient property of the closed fragment is that it is not sensitive to `branching'. Indeed, consider any Kripke frame
$\<W,<\>$ for $\mathsf{GL}^0$. To each $w\in W$ assign an ordinal $o(w)$ as follows: if $w$ is minimal, $o(w)=0$. Otherwise, $o(w)$ is the supremum of $o(v)+1$ over all $v<w$.

The map $o$ is well-defined because models of $\mathsf{GL}$ are well-founded. Further, because there are no variables, it is easy to check that $o:W\to\Lambda$ (where $\Lambda$ is a sufficiently large ordinal) is a bisimulation.

Thus to describe the modal logic of $W$ it is enough to describe $o(W)$. We can extend this idea to $\mathsf{GLP}_\Lambda$; if we have a well-founded frame $\mathfrak F=\langle W,\<<_\xi\>_{\xi<\Lambda}\rangle$, we can represent a world $w$ by the sequence $\vec o(w)=\<o_\xi(w)\>_{\xi<\Lambda}$, where $o_\xi$ is defined analogously to $o$. Thus we can identify elements of our model with sequences of ordinals. {It is a priori not clear that this representation suffices also for the polymodal case, and one of the main purposes of this paper is to see that it actually does.}


{Moreover}, there are certain conditions these sequences must satisfy. They arise from considering {\em worms}, which are formulas of the form $\<\xi_0\>...\<\xi_n\>\bot$. In various ways we can see worms as the backbone of the closed fragment of $\mathsf{GLP}$. It is known that each formula of ${\mathsf{GLP}^0_{\Lambda}}$ is equivalent to a Boolean combination of worms. Moreover, in \cite{paper0} it is shown that the axioms $\langle \alpha \rangle \phi \to \langle \beta \rangle \phi$ for $\alpha \geq \beta$ and $\langle \beta \rangle \phi \to [\alpha] \langle \beta \rangle \phi$ for $\alpha > \beta$ can be restricted simultaneously to worms to obtain an equivalent axiomatization of ${\mathsf{GLP}^0_{\Lambda}}$.

Given worms ${ A},{B}$ and an ordinal $\xi$, we define ${A\prec_{\xi}B}$ if $\vdash{B\to\<\xi\>A}$. This gives us a well-founded partial order. 

In \cite{WellOrders}, we study $\vec\Omega({A})=\langle\Omega_\xi(A)\rangle_{\xi<\Lambda}$, where $\Omega_\xi({A})$ is the order-type of $A$ under $<_\xi$. This gives us a good idea of what sequences may be included in the model; as it turns out, $\vec\Omega({A})$ is a `local bound' for $\vec o(w)$, giving rise to {\em $\le$-sequences} (see Section \ref{les}).

Our models naturally extend the model which was first defined and studied by Ignatiev for $\mathsf{GLP}^0_\omega$ in \cite{Ignatiev:1993}, which in our notation becomes $\mathfrak I^{\varepsilon_0}_\omega$, as well as Icard's topological variant, which here would be denoted $\mathfrak T^{\varepsilon_0}_\omega$. Originally, Ignatiev's study was an amalgamate of modal, arithmetical and syntactical methods. In \cite{Joosten:2004} the model was first submitted to a  purely modal analysis and \cite{BJV}  built forth on this work. In this paper, we prove soundness and completeness using purely semantic techniques, with the exception of a minor syntactic result from \cite{paper0} which is needed in the completeness proof.

\section{Hyperexponentials and -logarithms}\label{super}

In this section we shall introduce {\em hyperexponentials} and {\em hyperlogarithms} as a form of transfinite iteration of the function $-1+\omega^\xi$ {and its left-inverse $\le$, respectively}. These iterations have been used in \cite{WellOrders} for describing well-orders in the Japardize algebra, and will be essential in defining our semantics. We give only a very brief overview, but \cite{hyperexp} gives a thorough and detailed presentation.

We shall denote the class of all ordinals by $\mathsf{On}$ and the class of limit ordinals by $\mathsf{Lim}$.

\begin{definition}\label{hyperdef}
Let $\ex(\xi)=-1+\omega^\xi$. Then, we define the {\em hyperexponential} $\ex^\zeta\xi$ by the following recursion:
\begin{enumerate}
\item $\ex^0\xi=\xi$
\item $\ex^\xi 0=0$
\item $\ex^1=\ex$
\item $\ex^{\omega^\rho+\xi}=\ex^{\omega^\rho}\ex^\xi$, where $\xi<\omega^\rho+\xi$
\item $\ex^{\omega^\rho}(\xi+1)=\displaystyle\lim_{\zeta\to\omega^\rho}\ex^{\zeta}(\ex^{\omega^\rho}{(\xi)}+1)$, provided $\rho>0$
\item $\ex^{\omega^\rho}\xi=\displaystyle\lim_{\zeta\to\xi}\ex^{\omega^\rho}\zeta$ for $\xi\in\mathsf{Lim}$,  $\rho>0$.
\end{enumerate}
\end{definition}

\begin{prop}[Properties of {hyperexponentials}]\label{prophyp}
The family of functions $\langle\ex^\xi\rangle_{\xi\in\mathsf{On}}$ has the following properties:
\begin{enumerate}
\item $\ex^\xi$ is always a normal function\footnote{That is, strictly increasing and continuous.},
\item given ordinals $\xi,\zeta$, $\ex^{\xi+\zeta}=\ex^\xi \ex^\zeta$
\item given $\xi\in\mathsf{On}$, $\ex^{\xi+1}1=\displaystyle\lim_{n\to\omega}\ex^\xi n$
\item given $\xi\in\mathsf{Lim}$, $\ex^{\xi}1=\displaystyle\lim_{\zeta\to\xi}\ex^\zeta 1$
\item if $\xi<\zeta$ then $\ex^\xi\alpha\leq \ex^\zeta\alpha$
\item if $\xi+\zeta=\zeta$ then $\ex^\xi \ex^\zeta=\ex^\zeta$
\end{enumerate}
\end{prop}
\begin{definition}
For ordinals $\xi,\zeta$, define the {\em hyperlogarithms} $\le^\xi\zeta$ by the following recursion:
\begin{enumerate}
\item $\le^0\alpha=\alpha$
\item $\le^\xi 0=0$
\item $\le^\xi(\alpha+\beta)=\le^\xi\beta$ if $\beta>0$
\item $\le^{\omega^\rho+\xi}=\le^{\xi}\le^{\omega^\rho}$ provided $\xi<\omega^\rho+\xi$
\item $\le^{\omega^\rho}\ex^{\omega^\delta}\xi=
\begin{cases}
\le^{\omega^\rho}\xi&\text{if $\delta<\rho$}\\
\xi&\text{if $\delta=\rho$}\\
\ex^{\omega^\delta}\xi&\text{if $\delta>\rho$.}
\end{cases}$
\end{enumerate}
\end{definition}

Note in particular that, if $\xi=\zeta+\omega^\rho$, then $\ell\xi=\rho$; this is the {\em last exponent} or {\em end-logarithm} of $\xi$.

\begin{prop}[Properties of hyperlogarithms]
Hyperlogarithms have the following properties:
\begin{enumerate}
\item given ordinals $\xi,\zeta$, $\le^{\xi+\zeta}=\le^\zeta \le^\xi$
\item $\le^\zeta\alpha\leq \le^\xi\alpha$ whenever $\xi<\zeta$
\item $\le^{\zeta}\alpha=\le^\zeta \le^\xi\alpha$ whenever $\xi+\zeta=\zeta$.
\end{enumerate}
\end{prop}

Hyperlogarithms provide left-inverses for hyperexponentials:

\begin{lemm}\label{leftinv}
Given ordinals $\alpha,\beta,\xi$,
\begin{enumerate}
\item if $\alpha= \ex^\xi\beta$, then $\le^\xi\alpha=\beta$ and
\item if $\alpha<\ex^\xi\beta$, then $\le^\xi\alpha<\beta$.
\end{enumerate}

In general, if $\xi<\zeta$, then $\le^\xi\ex^\zeta=\ex^{-\xi+\zeta}$.
\end{lemm}

There is a close relation between the iterates $\ex^{\omega^\rho}\xi$ and Veblen functions; this is also described in detail in \cite{WellOrders}. For example, we have the following:

\begin{lemm}\label{rangelemm}
An ordinal $\xi$ lies in the range of $\ex^{\omega^\rho}$ if and only if, for all $\delta<\rho$, we have that $\xi=\ex^{\omega^\delta}\xi$. In particular, $\ex^{\omega^{\rho+1}}$ enumerates the fixpoints of $\ex^{\omega^{\rho}}$.
\end{lemm}

Like with Veblen functions, we may use hyperexponentials to give a sort of notation system for ordinals.

Given an ordinal $\xi$, say an expression
\[\xi=\sum_{i< I}\ex^{\alpha_i}\beta_i+n\]
is a {\em weak normal form} if $I,n<\omega$ and $0<\beta_i<\ex^{\alpha_i}\beta_i$ for all $i< I$. Note that weak normal forms are typically not unique; for example, we have
\[\omega^{\omega}=\ex^21=\ex\omega.\]

Say an ordinal $\xi$ is {\em definable} by a set $\Gamma$ if $\xi$ has a weak normal form
\[\sum_{i< I}\ex^{\alpha_i}\beta_i+n\]
where $n<\omega$, each $\alpha_i\in\Gamma$ and, inductively, $\Gamma$ defines each $\beta_i$. Every set of ordinals defines $0$.

Similar to Veblen normal forms, we have the following result:

\begin{prop}
Every ordinal $\xi$ has a weak normal form and hence is definable by $\Gamma$ large enough.
\end{prop}


\section{$\le$-sequences}\label{les}


In this section we shall describe the objects that are to be the `worlds' of our models. As stated before, these will be infinite sequences of ordinals; however, they must not only be weakly decreasing, but rather rapidly so. More specifically, they have to decrease at least as quickly as $\le^\xi$.

Given $\xi<\zeta$, we denote by $-\xi+\zeta$ the unique ordinal $\eta$ such that $\zeta=\xi+\eta$.

\begin{definition}[$\le$-sequence]
Let ${{\Theta}},\Lambda$ be ordinals.

We define a{n} {\em $\le$-sequence} (of depth ${{\Theta}}$ and length $\Lambda$) to be a function
\[f:\Lambda\to{{\Theta}}\]
such that, for every $\zeta\in(0,\Lambda)$, we have that \begin{equation}\label{dseq}
f(\zeta)\leq
\ell^{{-\xi+\zeta}} f(\xi)
\end{equation}
provided $\xi<\zeta$ is large enough.\footnote{
More precisely, given $\zeta\in(0,\Lambda)$ there is $\vartheta<\zeta$ such that (\ref{dseq}) holds whenever $\xi\in[\vartheta,\zeta)$.
}

If further
\[f(\xi+\zeta)=
\ell^{{\zeta}} f(\xi)\]
whenever $\xi+\zeta<\Lambda$ we say $f$ is {\em exact}.
\end{definition}

Let us see a few examples of $\le$-sequences:
\begin{itemize}
\item The sequence $f=\langle \omega^{\omega+1},\omega,1,0,\hdots\rangle$ is an $\le$-sequence, but it is not exact, since $\le\omega^{\omega+1}>\omega$. Note that once a sequence becomes zero it stabilizes, so we may represent sequences by their non-zero components.
\item The sequence $g=\langle \omega^{\omega+1},\omega+1,0,\hdots\rangle$ is an exact $\le$-sequence. Note that $g(1)>f(1)$ yet $g(2)<f(2)$.
\item The sequence $h$ given by
\[h(\xi)=
\begin{cases}
\varepsilon_0&\text{for $\xi<\omega$}\\
1&\text{for $\xi=\omega$}\\
0&\text{otherwise}
\end{cases}
\]
is an exact $\le$-sequence, since $\le^\omega\varepsilon_0=1$. Compare this to $h'$ defined as $h$ but with $h'(\omega)=0$; $h'$ is also an $\le$-sequence, but it is not exact.
\end{itemize}

As it turns out, to prove that an $\le$-sequence is exact, one only needs to check a fairly weak condition:

\begin{prop}\label{altd}
Let $f:\Lambda\to\Theta$. Then, the following are equivalent:
\begin{enumerate}
\item $f$ is exact;
\item for all $\zeta$ there is $\xi<\zeta$ such that
\[f(\zeta)=\le^{{-\xi+\zeta}}f(\xi).\]
\end{enumerate}
\end{prop}

\proof
A proof can be found in \cite{WellOrders}.
\endproof

Another nice property of exact sequences which will be useful later is the following:

\begin{lemm}\label{lastone}
If $f:\Lambda\to\Theta$ is an exact $\le$-sequence with $f(0)>0$ and $f(\xi)=0$ for some $\xi<\Lambda$, then there exists a maximum ordinal $\lambda$ such that $f(\lambda)\not=0$. Further, $f(\lambda)$ is a successor ordinal.
\end{lemm}

\proof
Let $\lambda$ be the supremum of all $\xi$ such that $f(\xi)>0$.

If $\lambda$ is a successor ordinal, then it immediately follows that $f(\lambda)>0$ (or $\lambda$ would not be the supremum).

Otherwise, write $\lambda=\gamma+\omega^\rho$ with $\rho>0$ and let $\vartheta\in[\gamma,\lambda)$ be large enough so that $f(\xi)=f(\vartheta)$ for all $\xi\in[\vartheta,\lambda)$; such a $\vartheta$ exists since $f$ is non-increasing.

Then, for $\delta<\rho$ we have that
\[\le^{\omega^\delta}f(\vartheta)=f(\vartheta+\omega^\delta)=f(\vartheta),\]
from which it follows that $f(\vartheta)=\ex^{\omega^\delta}f(\vartheta)$.

Hence $f(\vartheta)$ is a non-zero fixpoint of $\ex^{\omega^\delta}$ for all $\delta<\rho$, from which it follows using Lemma \ref{rangelemm} that it lies in the range of $\ex^{\omega^\rho}$ and thus is of the form $\ex^{\omega^\rho}\alpha$, with $\alpha>0$.

But then,
\[f(\lambda)=\le^{\omega^\rho}f(\vartheta)= \le^{\omega^\rho}\ex^{\omega^\rho}\alpha=\alpha\not=0.\]

Meanwhile, from maximality of $\lambda$ it follows that $f(\lambda+1)=\le f(\lambda)=0$, so $f(\lambda)$ must be a successor ordinal.
\endproof

We also have global characterizations for arbitrary $\le$-sequences:

\begin{prop}\label{tfae}
Given $f:\Lambda\to{{\Theta}}$, the following are equivalent:
\begin{enumerate}
\item $f$ is an $\le$-sequence
\item for every $\zeta\in(0,\Lambda)$,
\begin{enumerate}
\item if $\zeta=\xi+1$, $f(\zeta)\leq \le f(\xi)$ and
\item if $\zeta\in\mathsf{Lim}$, \[f(\zeta)\leq\lim_{\xi\to\zeta}\le^{{-\xi+\zeta}}f(\xi)\]
\end{enumerate}
\item for all $\xi<\zeta<\Lambda$,
\[\le f(\xi)\geq \le\ex^{{-\xi+\zeta}}f(\zeta)\]
\item for all $\xi<\zeta<\Lambda$,
\[\le f(\xi)\geq \le\ex^{\omega^{\le \zeta}}f(\zeta).\]
\end{enumerate}
\end{prop}

\proof
In principle 1 is stronger than 2; if (\ref{dseq}) holds for $\xi$ large enough, it holds in the limit. Note that for a successor ordinal $\zeta$, if any $\vartheta<\zeta$ exists such that (\ref{dseq}) holds for $\xi\in[\vartheta,\zeta)$, then we can always pick $\vartheta$ so that $\zeta=\vartheta+1$.

So we need only check that 2 implies 1 in the case of limit ordinals; but this follows from the fact that any function satisfying 2 must be non-increasing (which can be seen by a simple inspection) and thus limits are actually attained.

Likewise, 3 is in principle stronger than 4, because $\le e^{-\xi+\zeta} f(\zeta)\geq \le\ex^{\omega^{\le\zeta}}f(\zeta)$ independently of $\zeta,\xi$.

Thus our claim will be established if we show that 1 implies 3 and 4 implies 1.

Assume $f$ satisfies 1; let us check that it satisfies 3. For this we fix $\xi$ and proceed by induction on $\zeta$.

Write $\zeta=\gamma+\omega^\rho$ and pick $\vartheta<\zeta$ so that (\ref{dseq}) holds for all $\xi'\in[\vartheta,\zeta)$; without loss of generality, we can assume $\vartheta\geq \gamma$, so that $\zeta=\vartheta+\omega^\rho$. We may also assume that $\vartheta>\xi$, otherwise there is nothing to prove.

By induction on $\vartheta<\zeta$ we have that
\[\le f(\xi)\geq \le\ex^{-\xi+\vartheta} f(\vartheta).\]

Meanwhile, $f(\zeta)\leq \le^{\omega^\rho}f(\vartheta),$ so that by Lemma \ref{leftinv}.2, $\ex^{\omega^{\rho}} f(\zeta)\leq f(\vartheta).$

Thus
\[\le f(\xi)\geq \le\ex^{-\xi+\vartheta}\ex^{\omega^{\rho}} f(\zeta)=\le\ex^{-\xi+\vartheta+\omega^\rho} f(\zeta)=\le\ex^{-\xi+\zeta}f(\zeta),\]
which is what we wanted.

Finally, if $f$ satisfies 4, let us show that it also satisfies 1.

Choose $\zeta\in(0,\Lambda)$. Note that if $\zeta=\xi+1$ is a succesor, we can set $\vartheta=\xi$ and get $[\vartheta,\zeta)=\cbra\xi\cket$, while
\[\le f(\xi)\geq \le \ex^1 f(\zeta)=f(\zeta).\]
Thus we can assume otherwise and write $\zeta=\gamma+\omega^\rho$ with $\rho>0$. 

A quick inspection should show that $f$ is non-increasing, so we can pick $\vartheta\in[\gamma,\zeta)$ such that $f(\xi)=f(\vartheta)$ for all $\xi\in [\vartheta,\zeta)$. Because $\vartheta\geq\gamma$ we also have that ${-\xi+\zeta}=\omega^\rho$ for all such $\xi$. But by assumption
\[f(\xi)\geq{\le f(\xi)}\geq\le\ex^{\omega^{\rho}} f(\zeta)=\ex^{\omega^{\rho}} f(\zeta).\]
Hence our claim will follow if we show that $f(\xi)$ is in the range of $\ex^{\omega^{\rho}}$, since then we can apply $\le^{\omega^\rho}$ on both sides to obtain 
\[\le^{\omega^\rho}f(\xi)\geq f(\zeta).\]
To see this, pick $\delta<\rho$; in view of Lemma \ref{rangelemm} we must show that $f(\xi)$ is a fixpoint of $\ex^{\omega^{\delta}}$.

Since $f$ satisfies 4, we have that
\[{\le f(\xi)}\geq \le\ex^{\omega^{\delta}} f(\xi+\omega^\delta)=\ex^{\omega^{\delta}} f(\xi+\omega^\delta);\]
but $\delta<\rho$ so $\xi+\omega^\delta\in[\vartheta,\zeta)$, which implies that $f(\xi+\omega^\delta)=f(\xi)$ and thus this becomes ${\le f(\xi)}\geq \ex^{\omega^{\delta}} f(\xi)$. Hence $f(\xi)\geq \ex^{\omega^{\delta}} f(\xi);$
since $\ex^{\omega^{\delta}}$ is a normal function, it follows that $f(\xi)$ is a fixpoint of $\ex^{\omega^{\delta}}$, as claimed.
\endproof


\section{Generalized Ignatiev models}\label{gim}


Now rather than consider{ing} $\le$-se\-quen\-ces in isolation, we will be interested in forming a structure out of all $\le$-sequences {(possibly restricting depth and length)}. In this section we will generalize Ignatiev's universal model for $\glp_\omega^0$ to obtain models for $\glp_\Lambda^0$, independently of $\Lambda$. Our model combines ideas from Ignatiev's construction with results from previous sections to deal with limit modalities.

\begin{definition}[generalized Ignatiev model]
Given ordinals ${{\Theta}},\Lambda$, define a structure
\[\mathfrak I^{{\Theta}}_\Lambda=\<D^{{\Theta}}_\Lambda,\<<_\xi\>_{\xi<\Lambda}\>\]
by setting $D^{{\Theta}}_\Lambda$ to be the set of all $\le$-sequences of depth ${{\Theta}}$ and length $\Lambda$. Define $f<_\xi g$ if and only if $f(\zeta)=g(\zeta)$ for all $\zeta<\xi$ and $f(\xi)<g(\xi)$.
\end{definition}

Suppose $\Gamma$ is a set of ordinals and $\xi$ is any ordinal. We define the {\em $\Gamma$-norm} $\|\xi\|_{\Gamma}$ as the least $p<\omega$ such that one of the following holds:
\begin{enumerate}
\item $\xi=0$ and $p=0$,
\item $\xi=1$ and $p=1$,
\item $\xi=\alpha+\beta$ with $\alpha,\beta<\xi$ and $\|\alpha\|_{\Gamma}+\|\beta\|_{\Gamma}=p$ or
\item $\xi=\ex^{\gamma}\alpha$ with $\gamma\in\Gamma$ and $p=1+\|\alpha\|_{\Gamma}$.
\end{enumerate}

Let us compute a few examples:
\begin{itemize}
\item $2=1+1$ so $\|2\|_\Gamma=2$ independently of $\Gamma$. More generally, $\|n\|_\Gamma=n$ for $n<\omega$.
\item $\omega=\ex^1 1$ so $\|\omega\|_{\cbra 1\cket}=1+\|1\|_{\cbra 1\cket}=2$.
However, $\|\omega\|_{\varnothing}=\infty$, since $\omega$ cannot be written without the use of $\ex$.
\item $\|\varepsilon_0\|_{\cbra 1\cket}=\infty$, since $\varepsilon_0=\ex^{\omega} 1$ and $\varepsilon_0$ cannot be written with a smaller exponent.
\item $\|\omega^{\omega+1}\|_{\cbra 1\cket}=4$, since
\[
\begin{array}{lcl}
\|\omega^{\omega+1}\|_{\cbra 1\cket}&=&\|\ex^1(\omega+1)\|_{\cbra 1\cket}\\\\
&=&1+\|\omega+1\|_{\cbra 1\cket}\\\\
&=&1+\|\omega\|_{\cbra 1\cket}+\|1\|_{\cbra 1\cket}\\\\
&=&1+2+1\\\\
&=&4.
\end{array}\]
\end{itemize}

\begin{definition}[$\langle p,\Gamma\rangle$-approximation]\label{psa}
Given a natural number $p$ and a finite set of ordinals $\Gamma$, we say $\beta$ is a {\em $\langle p,\Gamma\rangle$-approximation} of $\alpha$ if $\beta<\alpha$ and $\|\beta\|_{\Gamma}\leq p$.
\end{definition}

Henceforth, we will say $\<p,\Gamma\>$ are {\em parameters} if $p<\omega$ and $\Gamma$ is a finite set of ordinals.

Clearly there are only finitely many $\langle p,\Gamma\rangle$-approximations of a given $\alpha$, and hence there is a maximum one: we denote it by $\lfloor\alpha\rfloor_{\Gamma}^p$. It will be convenient to stipulate $\lfloor 0\rfloor_{\Gamma}^p=-1$ for every $p,\Gamma$.

The approximations $\lfloor\alpha\rfloor_{\Gamma}^p$ will be very useful to us. One very elementary property they have is the following:

\begin{lemm}\label{unbounded}
If $\|\zeta\|_{\Gamma}\leq p$ and $\zeta<\xi$, then $\zeta\leq\lfloor \xi\rfloor^p_{\Gamma}$.
\end{lemm}

\proof
Obvious from Definition \ref{psa}.
\endproof

One can produce exact sequences from any function with finite domain, as we shall see.

Below, suppose that $r:\Gamma\to\Theta$, where $\Gamma\subseteq\Lambda$ is finite. We will define $\mathrm{dom}(r)$ to be the sequence $\vec\sigma$ which enumerates $\Gamma$. In general, if a function $s$ has domain $\Gamma$ we may write $s:\Gamma\to\Theta$ or $s:\vec\sigma\to\Theta$ indistinctly.

\begin{definition}[$\lceil r\rceil$]
Let $\vec\sigma=\langle\sigma_i\rangle_{i\leq I}$ be a finite, increasing sequence of ordinals containing zero with $\sigma_I<\Lambda$, $r:\vec\sigma\to \Theta$ be any function and $\delta_i=-\sigma_i+\sigma_{i+1}$.

Define a sequence $\lceil r\rceil:\Lambda\to\Theta$ by setting
\[
\lceil r\rceil(\xi)=
\begin{cases}
0&\text{for $\xi>\sigma_I$}\\\\
r(\sigma_I)+1&\text{for $\xi=\sigma_I$}\\\\
r(\sigma_i)+1+\ex^{\delta_{i}} \lceil r\rceil(\sigma_{i+1})&\text{for $\xi=\sigma_i$ with $i<I$}\\\\
\le^{\zeta}\lceil r\rceil(\sigma_{i})&\text{for $\xi=\sigma_i+\zeta<\sigma_{i+1}$.}
\end{cases}
\]
\end{definition}

Observe that this operation always produces exact $\le$-sequences:

\begin{lemm}\label{isad}
Given any finite $\Gamma\subseteq\Lambda$ and $r:\Gamma\to\Theta$, $\lceil r\rceil$ is an exact $\le$-sequence. 
\end{lemm}

\proof
We must establish that, given $\zeta<\Lambda$, there is $\xi<\zeta$ such that
\[\lceil r\rceil(\zeta)=\le^{-\xi+\zeta}\lceil r\rceil(\xi).\]

We make a few case distinctions:
\begin{description}
\item[$\lceil r\rceil(\zeta)=0$] Let $i$ be the largest index such that $\sigma_i<\zeta$ and take $\xi=\sigma_i$. Then, $\lceil r\rceil(\xi)$ is either zero or a successor ordinal and thus $\le^{-\xi+\zeta}\lceil r\rceil(\xi)=0=\lceil r\rceil(\zeta)$.

Note that this covers the case when $\xi>\sigma_I$.

\item[
$\zeta \in  (\sigma_i , \sigma_{i+1} {]}$
and $\lceil r\rceil(\zeta)>0$] 
Write $\zeta=\sigma_{i+1}+\zeta'$.

Then, $\lceil r\rceil(\zeta)=\le^{\zeta'} \lceil r\rceil(\sigma_i).$

\item[$\zeta=\sigma_{i+1}$] In this case,
\[\begin{array}{lcl}
\le^{-\sigma_i+\sigma_{i+1}}\lceil r\rceil(\sigma_{i}) &=&\le^{\delta_i}\lceil r\rceil(\sigma_{i})\\\\
&=&\le^{\delta_i}\left(r(\sigma_i)+1+\ex^{\delta_{i}} \lceil r\rceil(\sigma_{i+1})\right)\\\\
&=&\le^{\delta_i}\ex^{\delta_i}\lceil r\rceil(\sigma_{i+1})\\\\
&=&\lceil r\rceil(\sigma_{i+1}).\\
\end{array}\]
\end{description}\endproof

If $\vec\gamma=\langle \gamma_i\rangle_{i\leq I}$ is a finite, increasing sequence of ordinals with $\gamma_0=0$, we define $\Delta\vec\gamma=\cbra\delta_i\cket_{i<I}$, where $\delta_{i}=-\gamma_i+\gamma_{i+1}$.

\begin{lemm}\label{growth}
If $\vec\sigma=\langle\sigma_i\rangle_{i\leq I}$ is a finite, increasing sequence of ordinals below $\Lambda$, $r:\vec\sigma\to\Theta$ is any function and $\Gamma$ is a finite set of ordinals such that $\|r(\sigma_i)\|_\Gamma\leq p$ whenever $i\leq I$, then
\[\|\lceil r\rceil(\sigma_i) \|_{\Gamma\cup\Delta\vec\sigma}\leq (p+1)I\]
for all $i\leq I$.
\end{lemm}

\proof
One can show that
\[\|\lceil r\rceil(\sigma_i) \|_{\Gamma\cup\Delta\vec\sigma}\leq (p+1)(I-i)\]
by a simple backwards induction on $i$, observing the definition of $\lceil r\rceil(\sigma_i)$.
\endproof

The following simple, well-known lemma can be quite useful:

\begin{lemm}\label{less}
If $\alpha<\xi$ and $\beta\leq \le\xi$, then
\[\alpha+\omega^\beta\leq\xi.\]

\end{lemm}

\proof
By observation of the Cantor normal form of $\xi$.
\endproof

We can use constructions of the form $\lceil r\rceil$ to approximate $\le$-sequences. For this, the notion of a {\em radius} will be useful.

\begin{definition}
Given a function $f:\Lambda\to\Theta$, say another function $r:\Gamma\to \Theta$ is a {\em radius\footnote{The reason for this terminology should be clarified in Section \ref{topological}.} around $f$} if
\begin{enumerate}
\item $\Gamma\subseteq\Lambda$ is a finite set of ordinals containing $0$ and
\item if $r(\xi)$ is defined then $r(\xi)<f(\xi)$.
\end{enumerate}

\end{definition}
\begin{lemm}\label{between}
If $r$ is a function with finite domain $\vec\sigma=\langle\sigma_i\rangle_{i\leq I}$,
\begin{enumerate}
\item for every $i\leq I$, $r(\sigma_i)< \lceil r\rceil(\sigma_i)$ and
\item if $r$ is a radius around an $\le$-sequence $f$, then for all $\xi<\Lambda$, $\lceil r\rceil(\xi)\leq f(\xi).$
\end{enumerate}
\end{lemm}

\proof
Let $g=\lceil r\rceil$.

That $r(\sigma_i)< g(\sigma_i)$ is obvious from the definition of $\lceil r\rceil(\sigma_i)$, since it is always of the form
\begin{equation}\label{eqone}
r(\sigma)+\omega^\rho
\end{equation}
for some ordinal $\rho$.

To see the other inequality, we use backwards induction on $i$, noting that it is obvious when $\xi\geq\sigma_I$ or $g(\xi)=0$.

So we may suppose that $\xi\in[\sigma_i,\sigma_{i+1})$ and $g(\xi)>0$. Assume inductively that $f(\xi')\geq g(\xi')$ provided $\xi'\geq\sigma_{i+1}$.

By Proposition \ref{tfae}.3 we have that
\[
\le f(\xi)\geq \le \ex^{-\xi+\sigma_{i+1}} f(\sigma_{i+1})\geq\le \ex^{-\xi+\sigma_{i+1}}g(\sigma_{i+1})=\le g(\xi),
\]
where the second inequality follows by our induction hypothesis and the monotonicity of $\le\ex^{-\xi+\sigma_{i+1}}$.

Noting that $g(\xi)=\gamma+\omega^{\le g(\xi)}$ for some $\gamma<f(\xi)$ (possibly $\gamma=0$) we can then see using Lemma \ref{less} that
\[f(\xi)\geq \gamma+\omega^{\le g(\xi)}= g(\xi),\]
as claimed.
\endproof

There is another very natural operation to consider on $\le$-sequences, which under some conditions gives us new $\le$-sequences:

\begin{definition}[$\lambda$-concatenation]
Given sequences
\[f,g:\Lambda\to{{\Theta}},\]
we define their {\em $\lambda$-concatenation}
\[f\stackrel\lambda\ast g:\Lambda\to{{\Theta}}\]
by
\[
f\stackrel\lambda\ast g(\xi)=
\begin{cases}
f(\xi)&\text{if $\xi<\lambda$}\\
g(\xi)&\text{otherwise.}
\end{cases}
\]
\end{definition}

\begin{lemm}\label{concat}
If $f,g\in D^{{\Theta}}_\Lambda$ and $g(\lambda)\leq f(\lambda)$, then $f\stackrel\lambda\ast g$ is an $\le$-sequence.

If, further, $g(\lambda)< f(\lambda)$, then $f\stackrel\lambda\ast g<_\lambda f$.
\end{lemm}

\proof
Obvious from the definition of $f\stackrel\lambda\ast g$.
\endproof

We will often be interested in $r$ of a specific form. Given an $\le$-sequence $f$, a finite sequence of ordinals $\vec\sigma$ containing zero and $p<\omega$, define a radius $r=r[f,\vec\sigma,p]$ around $f$ with domain $\vec\sigma$ by
$r(\sigma_i)=\lfloor f(\sigma_i)\rfloor^p_{\Delta\vec\sigma}$. Then set
\[\lfloor f\rfloor^p_{\vec\sigma}=\lceil r[f,\vec\sigma,p]\rceil.\]

The sequence $\lfloor f\rfloor^p_{\vec\sigma}$ does not satisfy the same formulas of the modal language as $f$, but it does satisfy the same formulas that are `simple enough'. To see this we extend the notion of {\em $n$-bisimulation} to the slightly more general {\em $\langle p,\Gamma\rangle$-bisimulation}:

\begin{definition}[partial bisimulation]
Given $f,g\in D^{{\Theta}}_\Lambda$ and parameters $\langle p,\Gamma\rangle$, we say $f$ is {\em $\<p,\Gamma\>$-bisimilar} to $f$ (in symbols, $f\leftrightarroweq^p_{\Gamma}g$) by induction on $p$ as follows:

For $p=0$, any two $\le$-sequences are $\<p,\Gamma\>$-bisimilar.

For $p=q+1$, $f\leftrightarroweq^p_{\Gamma}g$ if and only if, for every $\gamma\in\Gamma$:
\begin{description}
\item[Forth] Whenever $f'<_\gamma f$, there is $g'<_\gamma g$ with $f'\leftrightarroweq^q_{\Gamma}g'$.
\item[Back] Whenever $g'<_\gamma g$, there is $f'<_\gamma f$ with $f'\leftrightarroweq^q_{\Gamma}g'$.
\end{description}
\end{definition}

The following is a well-known result from modal logic:

\begin{theorem}\label{bis}
If $\Gamma$ includes all modalities appearing in $\phi$ and $p$ is the modal depth of $\phi$, then whenever $\<\mathfrak I^\Theta_\Lambda,f\>\models\phi$ and $f\leftrightarroweq^p_{\Gamma} g$, it follows that $\<\mathfrak I^\Theta_\Lambda,g\>\models\phi$.
\end{theorem}

We may view $\Gamma$ indistinctly as a set or a sequence and thus also speak of $\langle p,\vec\sigma\rangle$-bisimulation. There is a close relation between $\langle q,\vec\sigma\rangle$-approximation and $\<p,\vec\sigma\>$-bisimulation, as we shall see.

Given $p<\omega$, a finite sequence of ordinals $\vec\sigma=\langle\sigma_i\rangle_{i\leq I}$ and $\le$-sequences $f,g$, say $f$ is {\em $\langle p,\vec\sigma\rangle$-close} to $g$, in symbols $f\sim^p_{\vec\sigma}g$, if for all $i\leq I$,
\[\lfloor f(\sigma_i)\rfloor^p_{\Delta\vec\sigma}=\lfloor g(\sigma_i)\rfloor^p_{\Delta\vec\sigma}.\]
It is not hard to check that $f\sim^p_{\vec\sigma}g$ if and only if $\lfloor g(\sigma_i)\rfloor^{p}_{\Delta\vec\sigma}<f(\sigma_i)$ and $\lfloor f(\sigma_i)\rfloor^{p}_{\Delta\vec\sigma}<g(\sigma_i)$ for all $i\leq I$. From this observation we easily obtain the following:

\begin{lemm}\label{issim}
Suppose that $f,g$ are $\le$-sequences, $p<\omega$ and $\vec\sigma=\langle\sigma_i\rangle_{i\leq I}$ a finite, increasing sequence of ordinals with $\sigma_0=0$.

Then, if for all $i\leq I$, either
\[\lfloor f(\sigma_i)\rfloor^p_{\Delta\vec\sigma}<g(\sigma_i)\leq f(\sigma_i)\]
or
\[\lfloor g(\sigma_i)\rfloor^p_{\Delta\vec\sigma}<f(\sigma_i)\leq g(\sigma_i),\]
it follows that $f\sim^p_{\vec\sigma}g$.
\end{lemm}

\proof
In the first case, we already have $\lfloor f(\sigma_i)\rfloor^p_{\Delta\vec\sigma}<g(\sigma_i),$ and from the inequality $g(\sigma_i)\leq f(\sigma_i)$ it immediately follows that $\lfloor g(\sigma_i)\rfloor^p_{\Delta\vec\sigma}<f(\sigma_i)$.

The second case is analogous, and since we obtained the desired inequalities for each $i\leq I$, we conclude that $f\sim^p_{\vec\sigma}g$, as claimed.
\endproof

\begin{lemm}\label{simulates}
Let $\vec\sigma$ be a finite sequence of ordinals. If $f,g$ are $\le$-sequences such that $f\sim^{(I+1)^p}_{\vec\sigma}g$, then $g\leftrightarroweq^p_{\vec\sigma}f$.
\end{lemm}

\proof
We prove the claim by induction on $p$. By symmetry it is enough to consider the `forth' condition.

Suppose that $f\sim^{(I+1)^{p+1}}_{\vec\sigma}g$; we will show that $f\leftrightarroweq^{p+1}_{\vec\sigma} g$.

Let $f'<_{\sigma_i}f$. We must find $g'<_{\sigma_i}g$ such that $f'\leftrightarroweq^p_{\vec\sigma}g'$; by induction hypothesis, it suffices to pick $g'$ such that $f'\sim^{(I+1)^p}_{\vec\sigma}g'$.

Let
\[g'=g\stackrel{\sigma_i}\ast \lfloor f'\rfloor^{(I+1)^p}_{\vec\sigma}.\]

First we must check that $g'$ is an $\le$-sequence and $g'<_{\sigma_i}f$. However, by Lemma \ref{concat}, it suffices to show that $g'(\sigma_i)<g(\sigma_i)$.

It follows from Lemma \ref{between} that $g'(\sigma_i)\leq f'(\sigma_i)$, and since $f'<_{\sigma_i}f$ we have that $g'(\sigma_i)<f(\sigma_i)$. But by Lemma \ref{growth},
\[\|\lfloor g'\rfloor^{(I+1)^p}_{\vec\sigma}(\sigma_i)\|_{\Delta\vec\sigma}\leq I(I+1)^p+I\leq (I+1)^{p+1},\]
so that $g'(\sigma_i)$ is an $\<(I+1)^{p+1},\Delta\vec\sigma\>$-approximation of $f(\sigma_i)$ and thus $g'(\sigma_i)\leq\lfloor f(\sigma_i)\rfloor^{(I+1)^{p+1}}_{\Delta\vec\sigma}$. Now, by assumption
\[\lfloor f(\sigma_i)\rfloor^{(I+1)^{p}}_{\Delta\vec\sigma}<g(\sigma_i),\]
so $g'(\sigma_i)<g(\sigma_i)$, as required.

We must also check that $f'\sim^{(I+1)^p}_{\vec\sigma}g'$; in other words, that for all $j\leq I$,
\[\lfloor f'(\sigma_j)\rfloor^{(I+1)^{p}}_{\Delta\vec\sigma}=\lfloor g'(\sigma_j)\rfloor^{(I+1)^{p}}_{\Delta\vec\sigma}.\]
But for $j<i$ this follows from the assumption that $g\sim^{(I+1)^{p+1}}_{\vec\sigma}f$, while for $j\geq i$ this follows form Lemmata \ref{between} and \ref{issim}.
\endproof

From this we immediately obtain the following:

\begin{cor}\label{notnew}
Any formula satisfiable over $\mathfrak I^\Theta_\Lambda$ is satisfied by an {\em exact} sequence $f\in D^\Theta_\Lambda$
\end{cor}

\proof
Suppose that $\<\mathfrak I^\Theta_\Lambda,g\>\models\phi$. Let $p$ be the modal depth of $\phi$ and $\vec\sigma$ a sequence of length $I$ which contains $0$ and every modality in $\phi$. Let $q=(I+1)^p$.

Then, by Lemmata \ref{between} and \ref{issim},
\[f=\lfloor g\rfloor^q_{\vec\sigma}\sim^q_{\vec\sigma}g,\]
so that by Lemma \ref{simulates}, $f\leftrightarroweq^p_{\vec\sigma}g$ and hence by Theorem \ref{bis}, $\<\mathfrak I^\Theta_\Lambda,f\>\models\phi$. Meanwhile, by Lemma \ref{isad}, $f$ is exact, as desired.
\endproof

\section{Generalized Icard topologies}\label{topological}

Corollary \ref{notnew} is a generalization of a known result; it has been observed in the past that Ignatiev's model has ``too many points'' in the sense that any formula can be satisfied on the {\em main axis}, i.e. the set of exact $\le$-sequences. However, these extra points are necessary if we wish to have Kripke semantics.

If we allow for topological semantics, then the main axis suffices.



We first note that the main axis of $\mathfrak I^{{\Theta}}_\Lambda$ can be identified with ${{\Theta}}$ in a canonical way, via the injection $\alpha\mapsto \vec\le\alpha$, where
\[\vec\le\alpha=\<\le^\xi\alpha\>_{\xi<\Lambda}.\]

Thus we can embed ${{\Theta}}$ into $D^{{\Theta}}_\Lambda$, and the image is precisely the main axis. Our goal for this section is to construct topologies $\mathcal T_\lambda$ for $\lambda<\Lambda$ which give us a polytopological model of $\glp^0_\Lambda$. For this, let us review the derived-set semantics of modal logic.

Recall that a {\em topological space} is a pair $\mathfrak X=\<X,\mathcal T\>$ where $\mathcal T\subseteq 2^X$ is a family of sets called `open' such that
\begin{enumerate}
\item $\varnothing,X\in\mathcal T$
\item if $U,V\in\mathcal T,$ then $U\cap V\in \mathcal T$ and
\item if $\mathcal U\subseteq\mathcal T$ then $\bigcup\mathcal U\subseteq\mathcal T$.
\end{enumerate}

Given $A\subseteq X$ and $x\in A$, we say $x$ is a {\em limit point} of $A$ if, given $U\in\mathcal T$ such that $x\in U$, we have that $(A\setminus \cbra x\cket)\cap U=\varnothing$. We denote the set of limit points of $A$ by ${d} A$, and call it the `derived set' of $A$.

We can define topological semantics for modal logic by interpreting Boolean operators in the usual way and setting
\[\lb\ps\psi\rb_\mathfrak X={d}\lb\psi\rb_\mathfrak X.\]

A {\em polytopological space} is a structure $\mathfrak X=\<X,\<\mathcal T_i\>_{i<I}\>,$ where each $\mathcal T_i$ is a topology. The derived set operator corresponding to $\mathcal T_i$ shall be denoted $d_i$. We can give conditions on the family of topologies so that $\glp_\Lambda$ is sound for $\mathfrak X$, and indeed $\glp_\omega$ is complete for these semantics \cite{BeklGabel:2011}.

Below, a topological space $\langle X,\mathcal T\rangle$ is {\em disperse} if every non-empty subset $A$ of $X$ has an isolated point; that is, if given $x\in A$ there is a neighborhood $U$ of $x$ (i.e., $x\in U\in\mathcal T$) such that $U\cap A=\cbra x\cket$.

We then have:

\begin{lemm}\label{glpspace}
Let $\Lambda$ be an ordinal and $\mathfrak X=\langle X,\<\mathcal T_\xi\>_{\xi<\Lambda}\rangle$ be a polytopological space.

Then,
\begin{enumerate}
\item L\"ob's axiom $[\xi]([\xi]\phi \to \phi)\to [\xi]\phi$ is valid on $\mathfrak X$ whenever $\<X,\mathcal T_\xi\>$ is disperse{,}

\item the axiom $\langle {\zeta} \rangle \phi \to \langle {\xi} \rangle \phi$ {for $\xi\leq\zeta$} is valid whenever $\mathcal T_\xi\subseteq \mathcal T_\zeta$ and

\item $\langle \xi \rangle \phi \to [\zeta] \langle \xi \rangle \phi$ {for $\xi<\zeta$} is valid if, whenever $A\subseteq X$, $d_\xi A\in \mathcal T_\zeta$.
\end{enumerate}

\end{lemm}
\proof
See, for example, \cite{BeklGabel:2011}.
\endproof

Although non-trivial spaces with these properties exist, they are hard to construct, and as in the case of Kripke semantics it turns out that restricting to the closed fragment significantly simplifies things.

Before defining our topological models, we recall the notion of a {\em subbasis}. Every collection of sets $\mathcal S\subseteq X$ such that $\bigcup \mathcal S= X$ gives rise to a least topology $\mathcal T$ containing $\mathcal S$. In this case we say $\mathcal S$ is a {\em subbasis} for $\mathcal T$. The elements of $\mathcal T$ are characterized as follows: $U\subseteq X$ is open if and only if, for every $x\in U$, there exists a {\em finite} subset $\mathcal N$ of $\mathcal S$ such that
\[x\in\bigcap \mathcal N\subseteq U.\] 

Finite intersections of sets in a subbasis are called {\em basic sets}\footnote{Of course not all bases are of this form, but this characterization will suffice for our purposes.}.

Our goal now is to build a sequence of topologies $\mathcal T_\lambda$ on ${{\Theta}}$ such that the resulting polytopological space is a model of $\glp_\Lambda$.

For this it will be convenient to assign three topological spaces to each ordinal $\xi$. We set:
\begin{enumerate}
\item $\xi^\bot$ to be $\xi$ with the trivial topology, i.e., the only opens are $\varnothing$ and all of $\xi$;
\item $\xi^{\mathcal I}$ to be $\xi$ with the initial segment topology, i.e. opens are intervals $[0,\gamma)$, with $\gamma\leq \xi+1$;
\item $\xi^{\mathcal O}$ to be $\xi$ with the {\em order topology}, i.e. with the topology generated by intervals of the form\footnote{Normally one defines order topologies using intervals of the form $(\alpha,\beta)$. But since we are dealing with ordinals, we can always rewrite $(\alpha,\beta]$ as $(\alpha,\beta+1)$, and thus intervals that are closed on the {\em right} are also open.} $(\alpha,\beta]$ or $[0,\beta]$ with $\beta\leq\xi+1$.
\end{enumerate}

For $\lambda<\Lambda$ define a topology $\mathcal T_\lambda$ on $|{{\Theta}}|^\Lambda$ (the bars indicate that exponentiation is taken set-theoretically, not as ordinals) by setting, for $\lambda<\Lambda$, $\mathcal T_\lambda$ to be the product topology
\[\prod_{\xi<\lambda}{{\Theta}}^{\mathcal O}\times{{\Theta}}^\mathcal I \times\prod_{\lambda<\xi}{{\Theta}}^\bot.\]

Note that $D^{{\Theta}}_\Lambda$ is a subset of $|{{\Theta}}|^\Lambda$ and, in turn, ${{\Theta}}$ can be seen as a subspace of $D^{{\Theta}}_\Lambda$ via the injection $\vec\le$. Hence $\mathcal T_\lambda$ induces a topology on ${{\Theta}}$ as a subspace of $|{{\Theta}}|^\Lambda$; we will not make a distinction and also denote this topology by $\mathcal T_\lambda$.

Equivalently, we can define $\mathcal T_\lambda$ by the subbasis consisting of intervals on coordinates below $\lambda$ and initial segments on $\lambda$. More precisely, subbasic sets are of the form
\[(\alpha,\beta]_\xi=\cbra f:\alpha< f(\xi)\leq\beta\cket\]
for some $\alpha<\beta\leq{{\Theta}}$ and $\xi<\lambda$, or of the form
\[[0,\beta]_\xi=\cbra f:f(\lambda)\leq\beta\cket\]
for $\xi\leq \lambda$.

We will call the resulting polytopological space $\mathfrak T^{{\Theta}}_\Lambda$. We interpret $\< \lambda\>$ by $\lb\<\lambda\>A\rb_{\mathfrak T^{{\Theta}}_\Lambda}={d}_\lambda\lb A\rb_{\mathfrak T^{{\Theta}}_\Lambda}$, i.e., the derivative operator with respect to the topology $\mathcal T_\lambda$.

There is a close connection between neighborhoods of $\xi$ and radii $r$ around $\vec\le\xi$. To see this, consider a $\mathcal T_\lambda$-neighborhood of $\xi$
\[U=\bigcap_{i\leq I}\left(\alpha_i,\le^{\sigma_i}\xi\right]_{\sigma_i}\cap[0,\le^\lambda\xi]_\lambda,\]
where all $\sigma_i<\lambda$; sets of this form form a basis for $\mathcal T_\lambda$.

Then define a radius $r$ around $f$ by $r(\sigma_i)=\alpha_i$. We can identify $U$ with $r$ and indeed will write $U=B^\lambda_r(\xi)$. Specifically, if $\mathrm{dom}(r)=\langle\sigma_0\rangle_{i\leq I}$ with $\sigma_I<\lambda$, set
\[B^\lambda_r(\xi)=\bigcap_{i\leq I}\left (r(\sigma_i),\le^{\sigma_i}\xi\right]_{\sigma_i}\cap [0,\le^\lambda\xi].\]
Thus we have a basis of $\mathcal T_\lambda$ such that neighborhoods of a point $\xi$ are identified with radii around $\vec\le\xi$.

Moreover, there is a sense in which $\mathcal T_\lambda$ is `irreflexive' in much the same way as $<_\lambda$:

\begin{lemm}\label{onlyyou}
Given $\xi\leq{{\Theta}}$ and $\lambda<\Lambda$, there is a $\mathcal T_\lambda$-neighborhood $U$ of $\xi$ such that whenever $\zeta\in U$ satisfies $\le^\lambda\zeta=\le^\lambda\xi$, it follows that $\zeta=\xi$.
\end{lemm}

\proof
By induction on $\lambda$.

If $\le^\lambda\xi=0$, let $\rho$ be the supremum of all $\zeta$ such that $\le^\zeta\xi>0$. By Lemma \ref{lastone}, $\rho$ is actually a maximum\footnote{Unless $\xi=0$, in which case the claim is trivial, as $\cbra 0\cket$ is open in all $\mathcal T_\lambda$.} and $\le^\rho\xi$ is a successor ordinal $\gamma+1$. Pick a $\mathcal T_\rho$-neighborhood $V$ of $\xi$ such that $\xi$ is the only element of $V$ with $\le^\rho\xi=\gamma+1$; such a neighborhood exists by induction hypothesis. It is not hard to check that $\xi$ is then the only element of $U=V\cap (\gamma,\gamma+1]_\rho\in \mathcal T_\lambda$, as desired.

Now assume that $\le^\lambda\xi>0$ and write $\lambda=\alpha+\omega^\rho$. By induction hypothesis, there is a $\mathcal T_\alpha$-neighborhood $V$ of $\xi$ such that there is no $\zeta\not=\xi$ in $V$ with $\le^\alpha\zeta= \le^\alpha\xi$.

Now write $\le^\alpha\xi=\gamma+\omega^{\beta}$, and let $U=V\cap (\gamma,\gamma+\omega^\beta]_\rho$; we claim that $U$ has the desired property.

Indeed, if $\zeta\in U$ has $\le^\lambda\zeta=\le^\lambda\xi$, this means that
\[\le^{\omega^\rho} \le^\alpha\zeta=\le^\lambda \xi,\]
i.e. $\le^\alpha\zeta$ is of the form $\delta+\ex^{\omega^{\rho}}\le^\lambda\xi$.

But we also know that
\[\le^{\omega^\rho} \le^\alpha\xi=\le^\lambda \xi,\]
so $\le^\alpha\xi$ is also of the form $\delta'+\ex^{\omega^\rho} \le^\lambda\xi$; in particular this implies that \[\omega^\beta=\ex^{\omega^{\rho}}\le^\lambda\xi.\]

But clearly the only element in the interval $(\gamma,\gamma+\omega^\beta]$ which is of this form is $\gamma+\omega^\beta$ itself, so it follows that
\[\le^\alpha\zeta=\gamma+\omega^\beta=\le^\alpha\xi.\]

By assumption $\xi$ was the only element in $U$ with this property, and we conclude that $\zeta=\xi$.
\endproof

We need one last simple definition before proving the main result of this section.

If $r,s$ are radii around $f$, define $t=r\sqcup s$ by
\[
t(\xi)=
\begin{cases}
r(\xi)&\text{if $r(\xi)$ is defined but $s(\xi)$ is not,}\\
s(\xi)&\text{if $s(\xi)$ is defined but $r(\xi)$ is not,}\\
\max\{ r(\xi),s(\xi)\}&\text{if $r(\xi)$ and $s(\xi)$ are both defined;}\\
\end{cases}
\]
everywhere else, $t(\xi)$ is undefined.

Then, one readily sees that
\[B^\lambda_{r\sqcup s}(f)=B^\lambda_{r}(f)\cap B^\lambda_{s}(f).\]

We are now ready to prove the following:
\begin{theorem}
Given $\xi<{{\Theta}}$ and a formula $\psi$,
\[\<\mathfrak T^{{\Theta}}_\Lambda,\xi\>\models \psi\Leftrightarrow\<\mathfrak I^{{\Theta}}_\Lambda,\vec \le\xi\>\models \psi.\]
\end{theorem}

\proof
We prove this by induction on $\psi$, where the cases for Booleans are trivial and we focus only on modal operators.

First assume that $\<\mathfrak I^{{\Theta}}_\Lambda,\vec \le\xi\>\models [\lambda]\psi.$ Let $\vec\sigma$ be an increasing sequence including $0$ as well as all modalities appearing in $\psi$ and let $p$ be the modal depth of $\psi$. Let $J$ be the largest index such that $\sigma_J<\lambda$.

Use Lemma \ref{onlyyou} to find a $\mathcal T_\lambda$-neighborhood $V$ of $\xi$ such that $\xi$ is the only element $\zeta$ in $V$ with $\le^\lambda\zeta=\le^\lambda\xi$, and let \[U=V\cap\bigcap_{i\leq J}\left(\lfloor \le^{\sigma_i}\xi\rfloor^{(I+1)^p}_{\Delta\vec\sigma},\le^{\sigma_i}\xi\right]_{\sigma_i}\cap [0,\le^\lambda\xi]_\lambda.\]

Let $\zeta\not=\xi\in U$ be arbitrary and consider  $f=\vec \le\xi\stackrel\lambda\ast \lfloor \vec\le\zeta\rfloor^{(I+1)^p}_{\vec\sigma}$. We know that $\zeta\in V$, so $f(\lambda)\leq\le^\lambda \zeta<\le^\lambda\xi$ and thus $f<_\lambda \vec \le\xi$; since we had assumed that $\<\mathfrak I^{{\Theta}}_\Lambda,\vec \le\xi\>\models[\lambda]\psi$, it follows that $\<\mathfrak I^{{\Theta}}_\Lambda,f\>\models\psi$.

For $i\leq J$ we obtain the inequality
\[\lfloor f(\sigma_i)\rfloor^{(I+1)^p}_{\Delta\vec\sigma}<\le^{\sigma_i}\zeta\leq f(\sigma_i)\]
from the assumption that $\zeta\in U$ and for $i>J$ we can see that
\[\lfloor \le^{\sigma_i}\zeta\rfloor^{(I+1)^p}_{\Delta\vec\sigma}<f(\sigma_i)\leq \le^{\sigma_i}\zeta\]
using Lemma \ref{between}; thus from Lemma \ref{issim} we have that $f\sim^{(I+1)^p}_{\vec\sigma}\vec\le\zeta$. It follows by Lemma \ref{simulates} and Theorem \ref{bis} that $\<\mathfrak I^{{\Theta}}_\Lambda,\vec\le\zeta\>\models\psi$ as well, and from our induction hypothesis, that $\<\mathfrak T^{{\Theta}}_\Lambda,\zeta\>\models\psi$. Since $\zeta\in U$ was arbitrary, we conclude that $\<\mathfrak T^{{\Theta}}_\Lambda,\xi\>\models[\lambda]\psi$, as claimed.

Now suppose that $\langle\mathfrak I^\Theta_\Lambda,\vec \le\xi\rangle\models\<\lambda\>\psi$, so that for some $g<_\lambda \vec \le\xi$ we have that $\<\mathfrak I^{{\Theta}}_\Lambda,g\>\models\psi$, and let $U$ be any $\mathcal T_\lambda$-neighborhood of $\xi$. Then, $U$ contains a neighborhood of $\xi$ of the form $B^\lambda_s(\xi)$ for some radius $s$ around $\vec\le\xi$.

Let $\vec\sigma=\langle\sigma_i\rangle_{i\leq I}$ be the sequence of all modalities in $\psi$ and $p$ be greater than the modal depth of $\psi$.

Let $r=r[f,\vec\sigma,p]$\footnote{I.e., for $i\leq I$, $r(\sigma_i)=\lfloor \le^{\sigma_i}\xi\rfloor^{(I+1)^p}_{\vec\sigma}$ and $r(\zeta)$ is undefined otherwise.} and $t=s\sqcup r$.

Now, let $h=\lceil t\rceil$. By Lemma \ref{isad}, $h$ is exact, so that $h=\vec\le \eta$ for $\eta=h(0)$. By Lemma \ref{between} we have that $\eta\in U$, since $h(\zeta)\in (s(\zeta),\le^\zeta\xi ]$ whenever $s(\zeta)$ is defined.

By Lemmata \ref{between} and \ref{issim} we have that $g\sim^{(I+1)^p}_{\vec\sigma}h$, since
\[h(\sigma_i)\in \left(\lfloor \le^\zeta\xi\rfloor^{(I+1)^p}_{\Delta\vec\sigma},\le^\zeta\xi\right]\]
for all $i\leq I$.

Therefore, by Lemma \ref{simulates}, $g\leftrightarroweq^p_{\vec\sigma}h$. It follows by Theorem \ref{bis} and induction on $\psi$ that $\<\mathfrak T^{{\Theta}}_\Lambda,\eta\>\models\psi$.

Since $U$ was arbitrary we conclude that
\[\<\mathfrak T^{{\Theta}}_\Lambda,\xi\>\models\<\lambda\>\psi.
\]
\endproof


\section{Soundness and completeness}\label{seccomp}

In this section we shall see that $\mathsf{GLP}^0_\Lambda$ is sound for both $\mathfrak I^\Theta_\Lambda$ and $\mathfrak T^\Theta_\Lambda$, as well as complete, provided that $\Theta$ is large enough.

Indeed, the soundness of the logics follows rather straightforwardly from our previous work.

\begin{theorem}
$\mathsf{GLP}^0_\Lambda$ is sound for both $\mathfrak I^\Theta_\Lambda$ and $\mathfrak T^\Theta_\Lambda$.
\end{theorem}

\proof
Most of the rules and axioms of $\mathsf{GLP}$ are standard, and we consider only the more unusual cases.

Note that, since $\mathfrak I^\Theta_\Lambda$ and $\mathfrak T^\Theta_\Lambda$ satisfy the same set of formulas, it suffices to check that each axiom is validated in {\em one} of these structures.

\begin{description}
\item[$\bra\xi\ket(\bra\xi\ket\phi\to\phi)\to\bra\xi\ket\phi$] This axiom is valid over $\mathfrak I^\Theta_\Lambda$ due to the transitivity and well-foundedness of $<_\xi$.
\item[$\<\zeta\>\phi\to\<\xi\>\phi$, provided $\xi<\zeta$] This axiom is valid over $\mathfrak T^\Theta_\Lambda$ because $\mathcal T_\zeta$ is a refinement of $\mathcal T_\xi$ (Lemma \ref{glpspace}.2).
\item[$\< \xi \> \phi\to\bra\zeta\ket \< \xi \> \phi$, provided $\xi<\zeta$] This axiom is valid over $\mathfrak I^\Theta_\Lambda$, given that if $f<_\zeta g$ and $f<_\xi h$, since $g(\vartheta)=f(\vartheta)$ for all $\vartheta<\zeta$, it follows that $h<_\xi g$.
\end{description}
\endproof

Before proceeding to consider completeness, let us see that `long' $\le$-sequences have large initial coordinates:

\begin{lemm}\label{almostthere}
Given ordinals $\lambda<\Lambda$ and $n<\omega$, there exists an $\le$-sequence $f$ with $f(\lambda)= n$ and $f(0)=\ex^\lambda n$; furthermore, if $g$ is any $\le$-sequence with $g(\lambda)\geq n$, then $g(0)\geq f(0)$.
\end{lemm}

\proof
First we shall construct an $\le$-sequence $f$ with $f(0)={\ex^\lambda}n$ and $f(\lambda)=n$, for any ordinal $\lambda$ and $n<\omega$.

Consider $f:\Lambda\to{{\Theta}}$ given by 
\[f(\xi)=
\le^{\xi}\ex^{\lambda}n.
\]
Clearly $f$ is an exact $\le$-sequence and, further,
\[f(\lambda)=\le^\lambda\ex^\lambda n=\ex^{-\lambda+\lambda}n=n,\]
while $f(0)=\ex^\lambda n$.

Now assume that $g(\lambda)\geq n$. By Proposition \ref{tfae}, $\le g(0)\geq \le\ex^\lambda n$, which clearly implies that $g(0)\geq \ex^\lambda n$.
\endproof

To continue we will need a syntactic result which is proven in \cite{paper0}. There it is stated in the more general setting of {\em Japaridze algebras}, which generalize both Kripke models and topological models; here we will state it for Kripke models, which are sufficient for our purposes.

\begin{theorem}\label{lastoneihope}
Let $\mathfrak F=\langle W,\langle R_\xi\rangle_{\xi<\Lambda}\rangle$ be a Kripke frame such that $\mathfrak F\models\glp^0_{\Lambda}$.

Then,
\begin{enumerate}
\item If $\Lambda=\lambda+1$ and for all $n<\omega$, $\langle\lambda\rangle^n\top$ is satisfied on $\mathfrak F$ or
\item if $\Lambda$ is a limit ordinal and for all $\lambda<\Lambda$, $\langle\lambda\rangle\top$ is satisfied on $\mathfrak F$,
\end{enumerate}
then $\glp^0_\Lambda$ is complete for $\mathfrak F$.
\end{theorem}

With this we may state and prove our main completeness result:

\begin{theorem}
The following are equivalent:
\begin{enumerate}
\item $\mathsf{GLP}^0_\Lambda$ is complete for $\mathfrak I^{{\Theta}}_\Lambda$
\item $\mathsf{GLP}^0_\Lambda$ is complete for $\mathfrak T^{{\Theta}}_\Lambda$
\item ${{\Theta}}\geq{\ex^\Lambda} 1.$
\end{enumerate}
\end{theorem}

\proof
Since $\mathfrak I^{{\Theta}}_\Lambda$ and $\mathfrak T^{{\Theta}}_\Lambda$ satisfy the same set of formulas, it suffices to show that 1 and 3 are equivalent.

First suppose that $\Lambda=\lambda+1$ is a successor ordinal.

Then, in view of Theorem \ref{lastoneihope}, $\glp^0_\Lambda$ is complete for $\mathfrak I^\Theta_\Lambda$ if and only if $\mathfrak I^\Theta_\Lambda$ satisfies $\<\lambda\>^n\top$ for all $n<\omega$. The latter is equivalent to the claim that, given $n<\omega$, there exists $f_n\in D^{{\Theta}}_\Lambda$ with $f_n(\lambda)\geq n$; by Lemma \ref{almostthere}.1, such an $f_n$ exists if and only if $\Theta>\ex^\lambda n$. But this must hold for all $n<\omega$, which by Proposition \ref{prophyp}.3 is equivalent to
\[{{\Theta}}\geq\lim_{n\to\infty}{\ex^\lambda}n =\ex^{\lambda+1} 1.\]

If $\Lambda$ is a limit ordinal, the argument is similar; here $\glp^0_\Lambda$ is complete for $\mathfrak I^\Theta_\Lambda$ if and only if it satisfies $\<\lambda\>\top$ for all $\lambda<\Lambda$, which is equivalent to the condition that, for all $\lambda<\Lambda$, there is $f_\lambda\in D^{{\Theta}}_\Lambda$ with $f_\lambda(\lambda)\geq 1$. By Lemma \ref{almostthere}.1, such an $f_\lambda$ exists if and only if $\Theta>{\ex^\lambda}1$. But this must hold for all $\lambda<\Lambda$; using Proposition \ref{prophyp}.4, this is equivalent to
\[{{\Theta}}\geq\lim_{\lambda\to\Lambda} {\ex^\lambda}1={\ex^\Lambda}1.\]
\endproof

\bibliographystyle{plain}
\bibliography{biblio}

\end{document}